\documentclass{article}
\usepackage{amsmath}

\newenvironment{proof}[1][Proof]{\noindent\textbf{#1.} }{\ \rule{0.5em}{0.5em}}
\input{tcilatex}

\begin{document}

\section{Introduction}

In [1] the author proved numerous results about ellipses inscribed\textbf{\ }%
in convex quadrilaterals. In particular, we proved that there exists a
unique ellipse of \textbf{minimal eccentricity}, $E_{I}$, inscribed in D.
This result applies to any convex quadrilateral, though the proof in [1]
assumes that D\ is not a trapezoid. In this paper, we discuss in detail the
special case of ellipses inscribed in parallelograms. In particular, in \S\ %
2 we give a direct proof(see Proposition 2) that there is a unique ellipse, $%
E_{I}$, of minimal eccentricity \textit{inscribed} in any given
parallelogram, D. Our main result in this regard is to give a geometric
characterization of $E_{I}$ for parallelograms(see Theorem 1), where we
prove that the smallest nonnegative angle between equal conjugate diameters
of $E_{I}$ equals the smallest nonnegative angle between the diagonals of D.
Similar results are known for the unique ellipse, $E_{A}$, of \textbf{%
maximal area} inscribed in a parallelogram, D(see, for example, [4]). Then
the equal conjugate diameters of $E_{A}$ are \textit{parallel} to the
diagonals of D. It is not too hard to prove this by proving the
corresponding result for the unit square and then using an affine
transformation. This works because of the affine invariance of the ratios of
corresponding areas. Since the eccentricity is not affine invariant, we
cannot reduce the problem of the minimal eccentricity ellipse inscribed in a
parallelogram to ellipses inscribed in squares.

In \S\ 3 we discuss ellipses inscribed in rectangles. We prove(see Theorem
2) that if $E_{M}$ is the unique ellipse inscribed in a rectangle, $R,$
which is tangent at the \textit{midpoints} of the sides of $R$, then $E_{M}$
is the unique ellipse of minimal eccentricity, maximal area, and maximal arc
length inscribed in $R$. While parts of Theorem 2 are known, this overall
characterization appears to be new. Of course, it then follows by affine
invariance that the unique ellipse of maximal area inscribed in a
parallelogram, D, is tangent at the midpoints of the sides of D. The other
parts of Theorem 2 do not hold in general for parallelograms, however.

In ([2], Proposition 1) the author proved that there is a unique ellipse, $%
E_{O}$, of minimal eccentricity circumscribed about any convex
quadrilateral, D. Also, in [2] the author defined D\ to be bielliptic if $%
E_{I}$ and $E_{O}$ have the same eccentricity. In \S\ 4 we show(Theorem 3)
that a parallelogram, D, is bielliptic if and only if the square of the
length of one of the diagonals of D\ equals twice the square of the length
of one of the sides of D.\qquad

Before proving our main results, we require the following lemma, which we
state without proof(see [6]).

\textbf{Lemma 1:} The equation $Ax^{2}+By^{2}+2Cxy+Dx+Ey+F=0,$ with $A,B>0,$
is the equation of an ellipse, $E_{0}$, if and only if $AB-C^{2}>0$ and $%
AE^{2}+BD^{2}+4FC^{2}-2CDE-4ABF>0$. Let $a$ and $b$ denote the lengths of
the semi--major and semi--minor axes, respectively, of $E_{0}$. Let $\phi $
denote the acute rotation angle of the axes of $E_{0}$ going
counterclockwise from the positive $x$ axis and let $\left(
x_{0},y_{0}\right) $ denote the center of $E_{0}$. Then 
\begin{equation}
a^{2}=\dfrac{AE^{2}+BD^{2}+4FC^{2}-2CDE-4ABF}{2(AB-C^{2})\left( A+B-\sqrt{%
(B-A)^{2}+4C^{2}}\right) },  \tag{1.1}
\end{equation}%
\begin{equation}
b^{2}=\dfrac{AE^{2}+BD^{2}+4FC^{2}-2CDE-4ABF}{2(AB-C^{2})\left( A+B+\sqrt{%
(B-A)^{2}+4C^{2}}\right) },  \tag{1.2}
\end{equation}%
and%
\begin{equation}
\phi =\dfrac{1}{2}\cot ^{-1}\left( \dfrac{A-B}{2C}\right) ,C\neq 0\text{ and 
}\phi =0\text{ if }C=0.  \tag{1.3}
\end{equation}

\section{Minimal Eccentricity}

\textbf{Lemma 2:} Let $Z$\ be the rectangle with vertices $%
(0,0),(l,0),(0,k), $and $(l,k),$ where $l,k>0$.

(A) The general equation of an ellipse, $\Psi $, inscribed in $Z$\ is given
by 
\begin{equation}
k^{2}x^{2}+l^{2}y^{2}-2l\left( k-2v\right)
xy-2lkvx-2l^{2}vy+l^{2}v^{2}=0,0<v<k.  \tag{2.1}
\end{equation}%
The corresponding points of tangency of $\Psi $ are 
\begin{equation}
\left( \dfrac{lv}{k},0\right) ,(0,v),\left( \dfrac{l}{k}\left( k-v\right)
,k\right) ,\text{and\ }(l,k-v).  \tag{2.2}
\end{equation}%
(B) If $a$ and $b$ denote the lengths of the semi--major and semi--minor
axes, respectively, of $\Psi $, then 
\begin{eqnarray}
a^{2} &=&\dfrac{2l^{2}\left( k-v\right) v}{k^{2}+l^{2}-\sqrt{%
(k^{2}+l^{2})^{2}-16l^{2}\left( k-v\right) v}}\text{ and}  \TCItag{2.3} \\
b^{2} &=&\dfrac{2l^{2}\left( k-v\right) v}{k^{2}+l^{2}+\sqrt{%
(k^{2}+l^{2})^{2}-16l^{2}\left( k-v\right) v}}.  \notag
\end{eqnarray}%
\textbf{Proof: }Let $S$ be the unit square with vertices $(0,0),(0,1),(1,0),$
and $(1,1)$. The map $T(x,y)=\left( \tfrac{1}{l}x,\tfrac{1}{k}y\right) $
maps $Z$ onto $S$ and $\Psi $ onto an ellipse, $T\left( \Psi \right) $.
Denote the points of tangency of $T\left( \Psi \right) $ with $S$ by $%
T_{1}=(t,0),T_{2}=(0,w),T_{3}=(s,1),$ and $T_{4}=(1,u),$ where $\left\{
t,w,s,u\right\} \subseteq (0,1)$. We may assume that the general equation of 
$T\left( \Psi \right) $ has the form $Ax^{2}+By^{2}+2Cxy+Dx+Ey+F=0$ with $%
A,B>0$. Since $T\left( \Psi \right) $ passes thru the points of tangency, we
have the equations 
\begin{eqnarray}
At^{2}+Dt+F &=&0,Bw^{2}+Ew+F=0  \TCItag{2.4} \\
As^{2}+B+2Cs+Ds+E+F &=&0,A+Bu^{2}+2Cu+D+Eu+F=0  \notag
\end{eqnarray}%
Using $y^{\prime }=-\tfrac{2Ax+2Cy+D}{2By+2Cx+E}$, $y^{\prime }\left(
T_{1}\right) =y^{\prime }\left( T_{3}\right) =0$ and the fact that the
tangents at $T_{2}$ and at $T_{4}$ are vertical, we also have the equations 
\begin{eqnarray}
2At+D &=&0,2Bw+E=0  \TCItag{2.5} \\
2As+2C+D &=&0,2Bu+2C+E=0  \notag
\end{eqnarray}%
Solving (2.4) and (2.5) for $B$ thru $F,s,t,$ and $u$ in terms of $A$ and $w$
yields $s=u=1-w,t=w,B=A,C=2Aw-A,D=-2Aw,E=-2Aw,F=Aw^{2}$. The equation of $%
T\left( \Psi \right) $ is then $x^{2}+y^{2}+2(2w-1)xy-2wx-2wy+w^{2}=0$. The
corresponding points of tangency of $T\left( \Psi \right) $ are $%
(w,0),(0,w),(1-w,1),(1,1-w)$. To obtain the corresponding equation of $\Psi ,
$ replace $x$ by $\tfrac{1}{l}x$ and $y$ by $\tfrac{1}{k}y$. That gives $%
\left( \tfrac{1}{l}x\right) ^{2}+\left( \tfrac{1}{k}y\right)
^{2}+2(2w-1)\left( \tfrac{1}{l}x\right) \left( \tfrac{1}{k}y\right)
-2w\left( \tfrac{1}{l}x\right) -2w\left( \tfrac{1}{k}y\right) +w^{2}=0$, or $%
k^{2}x^{2}+l^{2}y^{2}+2kl(2w-1)xy-2k^{2}lwx-\allowbreak
2kl^{2}wy+k^{2}l^{2}w^{2}=0$. The corresponding points of tangency of $\Psi $
are $T^{-1}(w,0)=\allowbreak \left( lw,0\right) $, $T^{-1}(0,w)=\allowbreak
\left( 0,kw\right) $, $T^{-1}(1-w,1)=\allowbreak \left( l\left( 1-w\right)
,k\right) $, and $T^{-1}(1,1-w)=\allowbreak \left( l,k\left( 1-w\right)
\right) $. Now let $v=kw$ to obtain (2.1) and (2.2). (2.3) follows easily
from Lemma 1, (1.1) and (1.2). 
\endproof%

We now prove a version of Lemma 2 for parallelograms.

\textbf{Proposition 1:} Let D\ be a parallelogram with vertices $%
O=(0,0),P=(l,0),Q=(d,k),$ and $R=(l+d,k)$, where $l,k>0,d\geq 0$.

(A) The general equation of an ellipse, $\Psi $, inscribed in D\ is given by%
\begin{gather}
k^{3}x^{2}+\left( k(d+l)^{2}-4dlv\right) y^{2}-2k\left( kd-2lv+kl\right) xy 
\tag{2.6} \\
-2k^{2}lvx+2klv\left( d-l\right) y+kl^{2}v^{2}=0,0<v<k.  \notag
\end{gather}%
(B) If $a$ and $b$ denote the lengths of the semi--major and semi--minor
axes, respectively, of $\Psi $, then 
\begin{equation}
\dfrac{b^{2}}{a^{2}}=1+\dfrac{m(v)+\allowbreak \left[ 4dlv-k\left(
(d+l)^{2}+k^{2}\right) \right] \sqrt{m(v)}}{8k^{2}l^{2}\left( k-v\right) v},
\tag{2.7}
\end{equation}%
where 
\begin{gather}
m(v)=\allowbreak 16l^{2}\left( d^{2}+k^{2}\right) v^{2}-8lk\left(
dk^{2}+d^{3}+2ld^{2}+l^{2}d+2k^{2}l\right) v+  \tag{2.8} \\
k^{2}\left( 2dl+l^{2}+d^{2}+k^{2}\right) ^{2}.  \notag
\end{gather}%
\textbf{Remark: }To be more precise, (A) means that any ellipse inscribed in
D\ has an equation of the form (2.6) for some $0<v<k$, and that any conic
with an equation of the form (2.6) for some $0<v<k$ defines an ellipse
inscribed in D.

\textbf{Proof: }Let $Z$\ be the rectangle with vertices $(0,0),(0,k),(l,0),$
and $(l,k)$. The map $T(x,y)=\left( x-\tfrac{d}{k}y,y\right) $ maps D\ onto $%
Z$. By Lemma 2, the general equation of $T\left( \Psi \right) $\ is given by
(2.1), with $x$ replaced by $x-\tfrac{d}{k}y$ and $y$ remaining the same.
That yields $k^{2}\left( x-\tfrac{d}{k}y\right) ^{2}+l^{2}y^{2}-2l\left(
k-2v\right) \left( x-\tfrac{d}{k}y\right) y-2lkv\left( x-\tfrac{d}{k}%
y\right) -2l^{2}vy+l^{2}v^{2}=0$, and some simplification gives (2.6). To
prove (B), by Lemma 1, (1.1) and (1.2), $\tfrac{b^{2}}{a^{2}}=\tfrac{(A+B)-%
\sqrt{(B-A)^{2}+4C^{2}}}{(A+B)+\sqrt{(B-A)^{2}+4C^{2}}}=\tfrac{\left[ (A+B)-%
\sqrt{(B-A)^{2}+4C^{2}}\right] ^{2}}{(A+B)^{2}-\left(
(B-A)^{2}+4C^{2}\right) }$, or%
\begin{equation}
\dfrac{b^{2}}{a^{2}}=\dfrac{(A+B)^{2}+(B-A)^{2}+4C^{2}-2(A+B)\sqrt{%
(B-A)^{2}+4C^{2}}}{\allowbreak 4(AB-C^{2})}.  \tag{2.9}
\end{equation}%
Let 
\begin{eqnarray}
A &=&k^{3},B=k(d+l)^{2}-4dlv,C=-k\left( kd-2lv+kl\right) ,  \TCItag{2.10} \\
D &=&-2k^{2}lv,E=2klv\left( d-l\right) ,\text{ and }F=kl^{2}v^{2}.  \notag
\end{eqnarray}%
Using (2.9) and (2.10), some simplification gives (2.7). 
\endproof%

\textbf{Proposition 2: }Let D\ be a parallelogram in the $xy$ plane. Then
there is a unique ellipse, $E_{I}$, of minimal eccentricity inscribed in D.

\textbf{Proof: }By using an isometry of the plane, we may assume that the
vertices of D\ are $O=(0,0),P=(l,0),Q=(d,k),$ and $R=(l+d,k)$, where $%
l,k>0,d\geq 0$. Let $E$ denote any ellipse inscribed in D\ and let $a$ and $b
$ denote the lengths of the semi--major and semi--minor axes, respectively,
of $E$. Let 
\begin{eqnarray}
g(v) &=&\dfrac{\allowbreak m(v)+\allowbreak \left[ 4dlv-k\left(
(d+l)^{2}+k^{2}\right) \right] \sqrt{m(v)}}{\left( k-v\right) v}, 
\TCItag{2.11} \\
h(v) &=&1+\dfrac{1}{8k^{2}l^{2}}g(v).  \notag
\end{eqnarray}%
By (2.7) of Proposition 1, $h(v)=\tfrac{b^{2}}{a^{2}}$. We shall now
minimize the eccentricity by maximizing\textbf{\ }$\tfrac{b^{2}}{a^{2}}$, or
equivalently by maximizing $g(v)$. Now $g^{\prime }(v)=0\iff $%
\begin{gather*}
\left( k-v\right) v\left[ m^{\prime }(v)+\allowbreak \left[ 4dlv-k\left(
(d+l)^{2}+k^{2}\right) \right] \dfrac{m^{\prime }(v)}{2\sqrt{m(v)}}+4dl\sqrt{%
m(v)}\right] - \\
\left[ m(v)+\allowbreak \left[ 4dlv-k\left( (d+l)^{2}+k^{2}\right) \right] 
\sqrt{m(v)}\right] (k-2v)=0
\end{gather*}%
$\iff $%
\begin{gather*}
\left( k-v\right) v\left[ 2\sqrt{m(v)}m^{\prime }(v)+\allowbreak \left[
4dlv-k\left( (d+l)^{2}+k^{2}\right) \right] m^{\prime }(v)+8dlm(v)\right] -
\\
2(k-2v)\sqrt{m(v)}\left[ m(v)+\allowbreak \left[ 4dlv-k\left(
(d+l)^{2}+k^{2}\right) \right] \sqrt{m(v)}\right] =0
\end{gather*}%
$\iff $%
\begin{gather*}
2\left( k-v\right) v\sqrt{m(v)}m^{\prime }(v)+\left( k-v\right) v\left( %
\left[ 4dlv-k\left( (d+l)^{2}+k^{2}\right) \right] m^{\prime
}(v)+8dlm(v)\right)  \\
-2(k-2v)\sqrt{m(v)}m(v)-2(k-2v)m(v)\allowbreak \left[ 4dlv-k\left(
(d+l)^{2}+k^{2}\right) \right] =0
\end{gather*}%
$\iff $%
\begin{gather*}
\left( k-v\right) v\left( \left[ 4dlv-k\left( (d+l)^{2}+k^{2}\right) \right]
m^{\prime }(v)+8dlm(v)\right) - \\
2(k-2v)m(v)\left( 4dlv-k\left( (d+l)^{2}+k^{2}\right) \right) = \\
\left[ 2(k-2v)m(v)-2v\left( k-v\right) m^{\prime }(v)\right] \sqrt{m(v)}
\end{gather*}%
$\iff $%
\begin{gather}
-\left[ 2\left( k^{2}+l^{2}+d^{2}\right) v-k\left(
2dl+l^{2}+d^{2}+k^{2}\right) \right] n(v)=  \tag{2.12} \\
\left[ 4dlv-k\left( (d+l)^{2}+k^{2}\right) \right] \left[ 2\left(
k^{2}+l^{2}+d^{2}\right) v-k\left( 2dl+l^{2}+d^{2}+k^{2}\right) \right] 
\sqrt{m(v)},  \notag
\end{gather}%
where $n(v)=m(v)\allowbreak +8k^{2}l^{2}\left( k-v\right) v$. If

\begin{equation*}
2\left( k^{2}+l^{2}+d^{2}\right) v-k\left( 2dl+l^{2}+d^{2}+k^{2}\right) \neq
0,
\end{equation*}%
then by (2.12), $g^{\prime }(v)=0\Rightarrow \sqrt{m(v)}=\tfrac{-n(v)}{%
4dlv-k\left( (d+l)^{2}+k^{2}\right) }\Rightarrow $

$\left( 4dlv-k\left( (d+l)^{2}+k^{2}\right) \right)
^{2}m(v)-n^{2}(v)=0\Rightarrow -64l^{4}v^{2}k^{4}\left( v-k\right)
^{2}=0\Rightarrow v=0$ or $v=k$. Since $0<v<k$ by assumption, that yields no
solution. Thus $g^{\prime }(v)=0$, and hence $h^{\prime }(v)=0$, if and only
if $2\left( k^{2}+l^{2}+d^{2}\right) v-k\left( 2dl+l^{2}+d^{2}+k^{2}\right)
=0\iff v=v_{\epsilon }$, where 
\begin{equation}
v_{\epsilon }=\dfrac{1}{2}k\dfrac{(d+l)^{2}+k^{2}}{k^{2}+d^{2}+l^{2}}. 
\tag{2.13}
\end{equation}%
It follows easily from l'Hospital's Rule that $\lim\limits_{v\rightarrow
0^{+}}g(v)=\lim\limits_{v\rightarrow k^{-}}g(v)=-8l^{2}k^{2}$, which implies
that $\lim\limits_{v\rightarrow 0^{+}}h(v)=\lim\limits_{v\rightarrow
k^{-}}h(v)=0$. Since $h(v)\geq 0$ for $0<v<k$, $h$ attains its' global
maximum at $v_{\epsilon }$ and the eccentricity is minimized when $%
v=v_{\epsilon }$. 
\endproof%

\textbf{Theorem 1:} Let $E_{I}$ denote the unique ellipse of minimal
eccentricity \textit{inscribed} in a parallelogram, D, in the $xy$ plane.
Then the smallest nonnegative angle between equal conjugate diameters of $%
E_{I}$ equals the smallest nonnegative angle between the diagonals of D.

\textbf{Proof: }As in the proof of Proposition 2, by using an isometry of
the plane, we may assume that the vertices of D\ are $%
O=(0,0),P=(l,0),Q=(d,k),$ and $R=(l+d,k)$, where $l,k>0,d\geq 0$. The
diagonals of D\ are $D_{1}=\overline{OR}$ and $D_{2}=\overline{PQ}$. We find
it convenient to define the following variables:%
\begin{eqnarray*}
G &=&(d+l)^{2}+k^{2},H=(d-l)^{2}+k^{2}, \\
J &=&d^{2}+k^{2}+l^{2},I=l^{2}-d^{2}-k^{2}.
\end{eqnarray*}%
There are three cases to consider: $I>0,I=0$(which implies that D\ is a
rhombus), and $I<0$. Assume first that $I>0$. Then $d^{2}+k^{2}<l^{2}$,
which implies that $d<l$ as well, and the lines containing $D_{1}$ and $D_{2}
$ have equations $y=\tfrac{k}{l+d}x$ and $y=\tfrac{k}{d-l}(x-l)$,
respectively. Let $\Psi $ equal the smallest nonnegative angle between $D_{1}
$ and $D_{2}$. We use the formula $\tan \Psi =\tfrac{m_{2}-m_{1}}{%
1+m_{1}m_{2}}$, where $m_{1}=\tfrac{k}{d-l}<m_{2}=\tfrac{k}{l+d}$. Some
simplification gives 
\begin{equation}
\tan \Psi =\dfrac{2kl}{I}.  \tag{2.14}
\end{equation}%
Let $E_{I}$ denote the the unique ellipse from Proposition 2 of minimal
eccentricity inscribed in D, and let $L$ and $L^{\prime }$ denote a pair of 
\textit{equal} conjugate diameters of $E_{I}$ . Let $a$ and $b$ denote the
lengths of the semi--major and semi--minor axes, respectively, of $E_{I}$.
It is known(see, for example, [5]) that $L$ and $L^{\prime }$ make equal
acute angles, on opposite sides, with the major axis of $E_{I}$. Let $\theta 
$ denote the acute angle going counterclockwise from the major axis of $E_{I}
$ to one of the equal conjugate diameters, which implies that $\tan \theta =%
\tfrac{b}{a}$. We shall show that $\tan ^{2}\left( 2\theta \right) =\tan
^{2}\Psi $, which will then easily yield $2\theta =\Psi $. By (2.7) of
Proposition 1, $\tfrac{b^{2}}{a^{2}}=h(v_{\epsilon })$, where $h(v)$ is
given by (2.11) and $v_{\epsilon }$ is given by (2.13). Thus 
\begin{equation*}
\tan \theta =\sqrt{h\left( v_{\epsilon }\right) }.
\end{equation*}%
By (2.13), 
\begin{gather*}
4dlv_{\epsilon }-k\left( (d+l)^{2}+k^{2}\right) =4dl\tfrac{1}{2}k\tfrac{%
2dl+l^{2}+d^{2}+k^{2}}{k^{2}+d^{2}+l^{2}}-k\left( (d+l)^{2}+k^{2}\right) = \\
-\tfrac{k\left( (d+l)^{2}+k^{2}\right) \left( (d-l)^{2}+k^{2}\right) }{%
d^{2}+k^{2}+l^{2}}=-\tfrac{kGH}{J},
\end{gather*}%
\begin{gather*}
v_{\epsilon }\left( k-v_{\epsilon }\right) =\tfrac{1}{2}k\tfrac{%
2dl+l^{2}+d^{2}+k^{2}}{k^{2}+d^{2}+l^{2}}\left( k-\tfrac{1}{2}k\tfrac{%
2dl+l^{2}+d^{2}+k^{2}}{k^{2}+d^{2}+l^{2}}\right) = \\
\tfrac{1}{4}k^{2}\left( 2dl+l^{2}+d^{2}+k^{2}\right) \tfrac{%
d^{2}-2dl+l^{2}+k^{2}}{\left( d^{2}+k^{2}+l^{2}\right) ^{2}}=\tfrac{1}{4}%
\allowbreak \tfrac{k^{2}\left( (d+l)^{2}+k^{2}\right) \left(
(d-l)^{2}+k^{2}\right) }{\left( d^{2}+k^{2}+l^{2}\right) ^{2}}=\tfrac{1}{4}%
\tfrac{k^{2}GH}{J^{2}},
\end{gather*}%
and by (2.8), after some simplification,%
\begin{gather*}
m\left( v_{\epsilon }\right) =m\left( \tfrac{1}{2}k\tfrac{%
2dl+l^{2}+d^{2}+k^{2}}{k^{2}+d^{2}+l^{2}}\right) = \\
\tfrac{k^{2}\left( (d+l)^{2}+k^{2}\right) \left( (d-l)^{2}+k^{2}\right)
\left( k^{2}-l^{2}+d^{2}\right) ^{2}}{\left( d^{2}+k^{2}+l^{2}\right) ^{2}}%
=\allowbreak \tfrac{k^{2}GHI^{2}}{J^{2}}.
\end{gather*}%
Hence $g(v_{\epsilon })=\left( \tfrac{k^{2}GHI^{2}}{J^{2}}-\tfrac{kGH}{J}%
\tfrac{k\sqrt{GH}I}{J}\right) \tfrac{4J^{2}}{k^{2}GH}=4I\left( I-\sqrt{G}%
\sqrt{H}\right) $, which implies, by (1.14), that $\allowbreak $%
\begin{equation}
h\left( v_{\epsilon }\right) =1+\dfrac{I\left( I-\sqrt{G}\sqrt{H}\right) }{%
2k^{2}l^{2}}  \tag{2.15}
\end{equation}%
Since $I^{2}-GH=\left( l^{2}-d^{2}-k^{2}\right) ^{2}-\left(
(d+l)^{2}+k^{2}\right) \left( (d-l)^{2}+k^{2}\right) =\allowbreak
-4l^{2}k^{2}<0,$ we have $I^{2}<GH$, which implies that $I-\sqrt{G}\sqrt{H}<0
$ and thus 
\begin{equation}
h\left( v_{\epsilon }\right) <1.  \tag{2.16}
\end{equation}%
Now $\tan 2\theta =\tfrac{2\tan \theta }{1-\tan ^{2}\theta }=\tfrac{2\sqrt{%
h\left( v_{\epsilon }\right) }}{1-h\left( v_{\epsilon }\right) }$, which
implies that 
\begin{gather}
\tan ^{2}2\theta =\tfrac{4h\left( v_{\epsilon }\right) }{\left( 1-h\left(
v_{\epsilon }\right) \right) ^{2}}=2\tfrac{2k^{2}l^{2}+I\left( I-\sqrt{G}%
\sqrt{H}\right) }{k^{2}l^{2}}\tfrac{4k^{4}l^{4}}{I^{2}\left( I-\sqrt{G}\sqrt{%
H}\right) ^{2}}=  \tag{2.17} \\
8k^{2}l^{2}\tfrac{2k^{2}l^{2}+I\left( I-\sqrt{G}\sqrt{H}\right) }{%
I^{2}\left( I-\sqrt{G}\sqrt{H}\right) ^{2}}.  \notag
\end{gather}

By (2.14) and (2.17), $\tan ^{2}2\theta =\tan ^{2}\Psi \iff 8k^{2}l^{2}%
\tfrac{2k^{2}l^{2}+I\left( I-\sqrt{G}\sqrt{H}\right) }{I^{2}\left( I-\sqrt{G}%
\sqrt{H}\right) ^{2}}=\tfrac{4k^{2}l^{2}}{I^{2}}\iff $

$4k^{2}l^{2}+2I\left( I-\sqrt{G}\sqrt{H}\right) =\left( I-\sqrt{G}\sqrt{H}%
\right) ^{2}\iff 4k^{2}l^{2}+2I^{2}-2I\sqrt{GH}=I^{2}-2I\sqrt{GH}+GH\iff $

$4k^{2}l^{2}+I^{2}=GH\iff 4k^{2}l^{2}+\left( l^{2}-d^{2}-k^{2}\right)
^{2}=\left( (d+l)^{2}+k^{2}\right) \left( (d-l)^{2}+k^{2}\right) $, which
holds for all $d,k,l\in \Re $. Thus $\tan ^{2}2\theta =\tan ^{2}\Psi $, and
since $\tan 2\theta >0$ and $\tan \Psi >0$ by (2.14) and (2.16), it follows
that $\tan 2\theta =\tan \Psi $.

\qquad Now suppose that $I=0$. One still has $d-l<0$, but now $\Psi =\tfrac{%
\pi }{2}$. One can let $I=0$ in (2.15) above by using a limiting argument.
Thus $h\left( v_{\epsilon }\right) =1$, which gives $2\theta =\tfrac{\pi }{2}
$. We omit the proof in the case when $I<0$. 
\endproof%

\textbf{Example: }Let $d=2$, $l=5$, and $k=4$, so that D\ is the
parallelogram with vertices $(0,0)$, $\left( 2,4\right) $, $\allowbreak
\left( 7,4\right) $, and $\allowbreak \left( 5,0\right) $. \ The minimal
eccentricity of ellipses inscribed in D\ is $\tfrac{65-\sqrt{65}}{65+\sqrt{65%
}}\approx \allowbreak 0.78$ and is attained with $v=\dfrac{50}{9}$. The
equation of $E_{I}$ is $1296x^{2}-531y^{2}+4464xy-18000x-13500\allowbreak
y+62500=0$. The common value of $2\theta $ and $\Psi $ equals $\tan
^{-1}8\approx \allowbreak 82.9^{\circ }$.

\textbf{Remark: }Theorem 1 does not extend in general to any convex
quadrilateral in the $xy$ plane. For example, consider the convex
quadrilateral with vertices $(0,0),(1,0),(0,1)$, and $(4,2)$. Using the
formulas from [1], one can show that there are two ellipses inscribed in D\
which satisfy $2\theta =\Psi $, but neither of those ellipses is the unique
ellipse of minimal eccentricity inscribed in D.

If D\ is a convex quadrilateral in the $xy$ plane, the line, $L$, thru the
midpoints of the diagonals of D\ plays an important role--it is the precise
locus of centers of ellipses inscribed in D. There is strong evidence that
the following is true.

\textbf{Conjecture: }Theorem 1 holds for any convex quadrilateral, D,\ with
the property that one of the diagonals of D\ is identical with $L$.

\qquad The details of a proof of this conjecture along the lines of the
proof of Theorem1 look messy. It is also possible that there is a similar
characterization for $E_{I}$ for any convex quadrilateral in the $xy$ plane.
Such a characterization would perhaps involve the angles between each
diagonal of D\ and between $L$ and each diagonal of D. However, we have not
found such a result which works with any examples.

\section{Rectangles}

The results in this paper have focused on ellipses of minimal eccentricity
inscribed in a parallelogram. We now discuss ellipses of minimal
eccentricity, maximal area, and maximal arc length inscribed in rectangles.
While some of the results in the following theorem are known, the overall
characterization is appears to be new.

\textbf{Theorem 2: }Let $Z$\ be a rectangle in the $xy$ plane. Then there is
a unique ellipse \textit{inscribed} in $Z$ which is tangent at the \textit{%
midpoints} of the four sides of $Z$, which we call the midpoint ellipse, $%
E_{M}$. $E_{M}$ has the following properties:

(A) $E_{M}$ is the unique ellipse of minimal eccentricity \textit{inscribed}
in $Z$.

(B) $E_{M}$ is the unique ellipse of maximal area \textit{inscribed} in $Z$.

(C) $E_{M}$ is the unique ellipse of maximal arc length \textit{inscribed}
in $Z$.

\textbf{Proof: }By using a translation, we may assume that the vertices of $Z
$\ are $O=(0,0),P=(l,0),Q=(0,k),$ and $R=(l,k)$, where $l,k>0$. Letting $v=%
\tfrac{1}{2}k$ in (2.2) shows the existence\textbf{\ }of an ellipse
inscribed in $Z$ which is tangent at the midpoints of the four sides of $Z$.
The fact that such an ellipse is unique\textbf{\ }follows easily and we omit
the proof. Now let $E$ denote any ellipse inscribed in $Z$\ and let $a$ and $%
b$ denote the lengths of the semi--major and semi--minor axes, respectively,
of $E$. To prove (A), as earlier we minimize the eccentricity by maximizing $%
\tfrac{b^{2}}{a^{2}}$. By (2.3), $\tfrac{b^{2}}{a^{2}}=\tfrac{k^{2}+l^{2}-%
\sqrt{(k^{2}+l^{2})^{2}-16l^{2}\left( k-v\right) v}}{k^{2}+l^{2}+\sqrt{%
(k^{2}+l^{2})^{2}-16l^{2}\left( k-v\right) v}}=8l^{2}h(v),0<v<k$, where 
\begin{equation*}
h(v)=\dfrac{v(k-v)}{\left( k^{2}+l^{2}\right) ^{2}+8l^{2}v(v-k)+\left(
k^{2}+l^{2}\right) \sqrt{\left( k^{2}+l^{2}\right) ^{2}+16l^{2}v(v-k)}}.
\end{equation*}%
A simple computation yields 
\begin{equation*}
h^{\prime }(v)=\tfrac{\allowbreak 8\left( k^{2}+l^{2}\right) l^{2}\left(
k-2v\right) }{\sqrt{\left( k^{2}+l^{2}\right) ^{2}+16l^{2}v(v-k)}\left(
\left( k^{2}+l^{2}\right) ^{2}+8l^{2}v\left( v-k\right) +(k^{2}+l^{2})\sqrt{%
\left( k^{2}+l^{2}\right) ^{2}+16l^{2}v(v-k)}\right) }.
\end{equation*}%
Thus $h^{\prime }(v)=0\iff v=\tfrac{1}{2}k$. Since $h(0)=h(k)=\allowbreak 0$
and $h(v)\geq 0$ for $0<v<k$, $h$ attains its' global maximum at $v=\tfrac{1%
}{2}k$ and the eccentricity is minimized when $E=E_{M}$. That proves (A). To
prove (B), we maximize the area of $E,\pi ab$, by maximizing $a^{2}b^{2}$.
By (2.3) again,%
\begin{gather*}
a^{2}b^{2}=-\tfrac{4\left( k-v\right) ^{2}l^{4}v^{2}}{\left( \sqrt{%
l^{4}+2l^{2}k^{2}+k^{4}+16l^{2}v^{2}-16l^{2}vk}-\left( k^{2}+l^{2}\right)
\right) \left( \sqrt{l^{4}+2l^{2}k^{2}+k^{4}+16l^{2}v^{2}-16l^{2}vk}+\left(
k^{2}+l^{2}\right) \right) } \\
=\tfrac{4\left( k-v\right) ^{2}l^{4}v^{2}}{16\left( k-v\right) l^{2}v}%
=\allowbreak \tfrac{l^{2}}{4}S(v),
\end{gather*}%
where $S(v)=\left( k-v\right) v$. It follows immediately that $S$ attains
its' global maximum at $v=\tfrac{1}{2}k$, which proves (B). To prove (C),
the arc length of $E$ is given by 
\begin{equation}
L=2\dint\limits_{0}^{\pi /2}\left[ a^{2}+b^{2}-\left( a^{2}-b^{2}\right)
\cos 2t\right] ^{1/2}dt.  \tag{3.1}
\end{equation}%
The proof we give is very similar to the proof in [3] that the ellipse of
maximal arc length\textbf{\ }inscribed in a square is a circle. Indeed, what
makes the proof work in [3] is that $a^{2}+b^{2}$ does not vary as $E$
varies over all ellipses inscribed in a square. For the rectangle, $Z$, 
\begin{gather*}
a^{2}+b^{2}=2l^{2}\left( k-v\right) v\left( \tfrac{1}{k^{2}+l^{2}-\sqrt{%
(k^{2}+l^{2})^{2}-16l^{2}\left( k-v\right) v}}+\tfrac{1}{k^{2}+l^{2}+\sqrt{%
(k^{2}+l^{2})^{2}-16l^{2}\left( k-v\right) v}}\right)  \\
=2l^{2}\left( k-v\right) v\tfrac{2\left( k^{2}+l^{2}\right) }{16l^{2}v\left(
k-v\right) }=\tfrac{1}{4}\left( k^{2}+l^{2}\right) ,
\end{gather*}%
which of course does not vary as $E$ varies over all ellipses inscribed in $Z
$. Now 
\begin{gather*}
a^{2}-b^{2}=\allowbreak \tfrac{4\left( k-v\right) l^{2}v\sqrt{%
l^{4}+2l^{2}k^{2}+k^{4}+16l^{2}v^{2}-16l^{2}vk}}{\left( k^{2}+l^{2}-\sqrt{%
l^{4}+2l^{2}k^{2}+k^{4}+16l^{2}v^{2}-16l^{2}vk}\right) \left( k^{2}+l^{2}+%
\sqrt{l^{4}+2l^{2}k^{2}+k^{4}+16l^{2}v^{2}-16l^{2}vk}\right) } \\
=4\left( k-v\right) l^{2}v\tfrac{\sqrt{%
l^{4}+2l^{2}k^{2}+k^{4}+16l^{2}v^{2}-16l^{2}vk}}{\allowbreak 16\left(
k-v\right) l^{2}v}\allowbreak =\tfrac{1}{4}\sqrt{g(v)},
\end{gather*}%
where $\allowbreak g(v)=\left( k^{2}+l^{2}\right) ^{2}-16l^{2}v\left(
k-v\right) $. Hence by (3.1), $L=L(v)=\tfrac{1}{2}\dint\limits_{0}^{\pi /2}%
\left[ k^{2}+l^{2}-\sqrt{g(v)}\cos 2t\right] ^{1/2}dt$. As in [3], splitting
the integral up and making a change of variable gives 
\begin{equation*}
L(v)=\dint\limits_{0}^{\pi /4}\left[ \left( k^{2}+l^{2}-\sqrt{g(v)}\cos
2t\right) +\left( k^{2}+l^{2}+\sqrt{g(v)}\cos 2t\right) \right] ^{1/2}dt
\end{equation*}%
Let $p=k^{2}+l^{2}$ and $u(v,t)=\sqrt{g(v)}\cos 2t$, which gives 
\begin{equation}
L(v)=\dint\limits_{0}^{\pi /4}\left[ (p-u(v,t))^{1/2}+(p+u(v,t))^{1/2}\right]
dt.  \tag{3.2}
\end{equation}%
Now $g$ attains its global minimum on $(0,k)$ when $v=\tfrac{1}{2}k$. Thus,
for each $0<t<\tfrac{\pi }{4}$, $u(v,t)\geq u\left( \tfrac{1}{2}k,t\right) $%
, with equality if and only if $v=\tfrac{1}{2}k$. Also, the function $%
f(x)=(p-x)^{1/2}+(p+x)^{1/2}$ is strictly decreasing for $0<x<p$(see [3]).
Hence, for each $0<t<\tfrac{\pi }{4}$, $(p-u(v,t))^{1/2}+(p+u(v,t))^{1/2}%
\leq (p-u\left( \tfrac{1}{2}k,t\right) )^{1/2}+(p+u\left( \tfrac{1}{2}%
k,t\right) )^{1/2}$, again with equality if and only if $v=\tfrac{1}{2}k$.
Thus by (3.2), $L(v)$ attains its' unique maximum on $(0,k)$ when $v=\tfrac{1%
}{2}k$. 
\endproof%

\textbf{Remark: }There is no ellipse of minimal arc length\textbf{\ }%
inscribed in $Z$.

\textbf{Remark: }Showing that there is a unique ellipse of maximal arc length%
\textbf{\ }inscribed in a general convex quadrilateral and/or characterizing
such an ellipse appears to be a very nontrivial problem. Even for
parallelograms it appears to be difficult. In general, $a^{2}+b^{2}$ does
not remain constant as $E$ varies over all ellipses inscribed in a
parallelogram. Numerical evidence suggests strongly that the ellipse of of
minimal eccentricity inscribed in a parallelogram, D, is \textbf{not} the
ellipse of maximal arc length inscribed in D.

\section{Bielliptic Parallelograms}

Let D\ be a convex quadrilateral. In ([1], Theorem 4.4) the author proved
that there is a unique ellipse, $E_{I}$, of minimal eccentricity \textit{%
inscribed} in D. In ([2], Proposition 1) we also proved that there is a
unique ellipse, $E_{O}$, of minimal eccentricity \textit{circumscribed}
about D. In [2] the author defined D\ to be bielliptic if $E_{I}$ and $E_{O}$
have the \textbf{same eccentricity}. This generalizes the notion of
bicentric quadrilaterals, which are quadrilaterals which have both a
circumscribed and an inscribed circle. In [2] we gave an example of a
bielliptic convex quadrilateral which is not a parallelogram and which is
not bicentric. Of course every square is bicentric. For parallelograms in
general we prove the following.

\textbf{Theorem 3: }A parallelogram, D, is bielliptic if and only if the
square of the length of one of the diagonals of D\ equals twice the square
of the length of one of the sides of D.

\textbf{Proof: }We prove the case when D\ is \textbf{not} a rectangle, in
which case the proof below can be modified to show that D\ is bielliptic if
and only if it's a square, which certainly satisfies the conclusion of
Theorem 3. Then, by using an isometry of the plane, we may assume that the
vertices of D\ are $O=(0,0),P=(l,0),Q=(d,k),$ and $R=(l+d,k)$, where $d,k,l>0
$. It is not hard to show that 
\begin{equation}
kux^{2}+ky^{2}-2udxy-klux+\left[ ud(l+d)-k^{2}\right] y=0,0<u<\dfrac{k^{2}}{%
d^{2}}  \tag{4.1}
\end{equation}%
is the general equation of an ellipse passing thru the vertices of D. We
leave the details to the reader. By Lemma 1, it follows that $\tfrac{b^{2}}{%
a^{2}}=\tfrac{k(u+1)-\sqrt{k^{2}(1-u)^{2}+4d^{2}u^{2}}}{k(u+1)+\sqrt{%
k^{2}(1-u)^{2}+4d^{2}u^{2}}}=h(u)$, where 
\begin{equation*}
h(u)=\tfrac{\left( k(u+1)-\sqrt{k^{2}(1-u)^{2}+4d^{2}u^{2}}\right) ^{2}}{%
4u\left( k^{2}-ud^{2}\right) }.
\end{equation*}%
Differentiating with respect to $u$, it follows that $h^{\prime }(u)=0,0<u<%
\tfrac{k^{2}}{d^{2}},$ if and only if $u=\tfrac{k^{2}}{k^{2}+2d^{2}}$.
Substituting yields $h\left( \tfrac{k^{2}}{k^{2}+2d^{2}}\right) =\tfrac{%
\left( d^{2}+k^{2}-d\sqrt{d^{2}+k^{2}}\right) ^{2}}{k^{2}\left(
d^{2}+k^{2}\right) }$, and simplifying gives 
\begin{equation}
1-\dfrac{b^{2}}{a^{2}}=2d\dfrac{\sqrt{d^{2}+k^{2}}-d}{k^{2}}  \tag{4.2}
\end{equation}%
for the unique ellipse of minimal eccentricity, $E_{O}$, circumscribed about
D. As in the proof of Theorem 1, there are three cases to consider: $%
I>0,I=0,I<0$, where $I=l^{2}-d^{2}-k^{2}$. Assume first that $I>0$. Then by
(2.15) in the proof of Theorem 1, $\tfrac{b^{2}}{a^{2}}=h(v_{\epsilon })=1+%
\tfrac{I\left( I-\sqrt{G}\sqrt{H}\right) }{2k^{2}l^{2}}$ for the unique
ellipse, $E_{I}$, of minimal eccentricity inscribed in D, where $%
G=(d+l)^{2}+k^{2}$ and $H=(d-l)^{2}+k^{2}$. Setting the eccentricities of $%
E_{I}$ and $E_{O}$ equal is equivalent to 
\begin{equation}
2d\dfrac{\sqrt{d^{2}+k^{2}}-d}{k^{2}}=\dfrac{I\left( \sqrt{G}\sqrt{H}%
-I\right) }{2k^{2}l^{2}}.  \tag{4.3}
\end{equation}%
Then (4.3) holds if and only if $4l^{2}d\left( \sqrt{d^{2}+k^{2}}-d\right)
=I\left( \sqrt{G}\sqrt{H}-I\right) \iff \sqrt{G}\sqrt{H}I=\allowbreak
4l^{2}d\left( \sqrt{d^{2}+k^{2}}-d\right) +I^{2}\iff $

$GH\allowbreak I^{2}-\left( \allowbreak 4l^{2}d\left( \sqrt{d^{2}+k^{2}}%
-d\right) +I^{2}\right) ^{2}\allowbreak =0\iff 4l^{2}\left( G-2l^{2}\right)
\left( H-2l^{2}\right) \allowbreak \left( 2d^{2}+k^{2}-2d\sqrt{d^{2}+k^{2}}%
\right) =0\iff $ one of the following equations holds:

\begin{eqnarray}
k^{2}+d^{2}-2dl-l^{2} &=&0  \TCItag{4.4} \\
k^{2}+d^{2}+2dl-l^{2} &=&0  \TCItag{4.5} \\
2d^{2}+k^{2}-2d\sqrt{d^{2}+k^{2}} &=&0  \TCItag{4.6}
\end{eqnarray}%
But $I>0\Rightarrow d^{2}+k^{2}-l^{2}<0\Rightarrow $ (4.4) cannot hold. $%
\left( 2d^{2}+k^{2}\right) =2d\sqrt{d^{2}+k^{2}}\Rightarrow \left(
2d^{2}+k^{2}\right) ^{2}-4d^{2}\left( d^{2}+k^{2}\right) =0\Rightarrow
k^{4}=0\Rightarrow k=0$, and thus (4.6) cannot hold either. We are left with
(4.5), which is a valid equation. Now if $I<0$, one then obtains, exactly
along the same lines, (4.4). Thus the eccentricities of $E_{I}$ and $E_{O}$
are equal if and only if (4.4) or (4.5) holds.

\qquad The diagonals of D\ are $D_{1}=\overline{OR}$ and $D_{2}=\overline{PQ}
$, and thus the squares of the lengths of one of the diagonals are $%
\left\vert D_{1}\right\vert ^{2}=(l+d)^{2}+k^{2}$ and $\left\vert
D_{2}\right\vert ^{2}=(l-d)^{2}+k^{2}$. The squares of the lengths of the
sides are $\left\vert \overline{OQ}\right\vert ^{2}=d^{2}+k^{2}$ and $%
\left\vert \overline{OP}\right\vert ^{2}=l^{2}$. Now $\left\vert
D_{1}\right\vert ^{2}=2\left\vert \overline{OQ}\right\vert ^{2}\iff
(l+d)^{2}+k^{2}=2d^{2}+2k^{2}\iff k^{2}+d^{2}-2dl-l^{2}=0$, which is (4.4).
Similarly, $\left\vert D_{1}\right\vert ^{2}=2\left\vert \overline{OP}%
\right\vert ^{2}\iff (l+d)^{2}+k^{2}=2l^{2}\iff d^{2}+2dl-l^{2}+k^{2}=0$,
which is (4.5). One can easily check that $\left\vert D_{2}\right\vert
^{2}=2\left\vert \overline{OQ}\right\vert ^{2}$ or $\left\vert
D_{2}\right\vert ^{2}=2\left\vert \overline{OQ}\right\vert ^{2}$ yields
(4.4) or (4.5) as well.

\qquad Finally suppose that $I=0$. Letting $I$ approach $0$ in (1.18) shows
that the unique ellipse, $E_{I}$, of minimal eccentricity inscribed in D\ is
a circle, which has eccentricity $0$. But $1-\dfrac{b^{2}}{a^{2}}=0$ in
(4.2) if and only if $d=0$. In that case D\ is a square, which again
satisfies the conclusion of Theorem 3. 
\endproof%

\textbf{Remark: }In the proof above, $\left\vert D_{1}\right\vert
^{2}+\allowbreak \left\vert D_{2}\right\vert ^{2}=\allowbreak 2\left(
d^{2}+l^{2}+k^{2}\right) =2\left\vert \overline{OP}\right\vert
^{2}+2\left\vert \overline{OQ}\right\vert ^{2}$ for \textit{any }%
parallelogram, D, and not just a bielliptic parallelogram. Hence if, say, $%
\left\vert D_{1}\right\vert ^{2}=2\left\vert \overline{OP}\right\vert ^{2}$,
then it follows automatically that $\allowbreak \left\vert D_{2}\right\vert
^{2}=2\left\vert \overline{OQ}\right\vert ^{2}$.

\textbf{Example: }Let $l=6$, $k=2\sqrt{2}$, and $d=2$. Then $I=\allowbreak
24>0$, and the common eccentricity of $E_{I}$ and $E_{O}$ is $\sqrt{3}-1$.
The squares of the lengths of one of the diagonals are $\allowbreak 72$ and $%
\allowbreak 24$, and the squares of the lengths of the sides are $%
\allowbreak 12$ and $\allowbreak 36$.

$v=\allowbreak \tfrac{3}{2}\sqrt{2}$ yields the ellipse, $E_{I}$, of minimal
eccentricity inscribed in D, and $4\sqrt{2}x^{2}+14\sqrt{2}y^{2}+4xy-36\sqrt{%
2}x-72y+81\sqrt{2}\allowbreak =0$ is the equation of $E_{I}\allowbreak $. $%
u=\allowbreak \tfrac{1}{2}$ yields the ellipse, $E_{O}$, of minimal
eccentricity circumscribed about D, and $\sqrt{2}x^{2}+2\sqrt{2}y^{2}-2xy-6%
\sqrt{2}x=0$ is the equation of $E_{O}$.

\section{References}

[1] Alan Horwitz, \textquotedblleft Ellipses of maximal area and of minimal
eccentricity inscribed in a convex quadrilateral\textquotedblright ,
Australian Journal of Mathematical Analysis and Applications, 2(2005), 1-12.

[2] Alan Horwitz, \textquotedblleft Ellipses of minimal area and of minimal
eccentricity circumscribed about a convex quadrilateral\textquotedblright ,
Australian Journal of Mathematical Analysis and Applications, to appear.

[3] The William Lowell Putnam Mathematical Competition problems and
solutions, Mathematical Association of America, 1965-1984.

[4] \texttt{Online }http://www.math.uoc.gr/\symbol{126}%
pamfilos/eGallery/problems/ParaInscribedEllipse.html

[5] George Salmon, A treatise on conic sections, 6th edition, Chelsea
Publishing Company, New York.

[6] \texttt{Online http://mathworld.wolfram.com/Ellipse.html}

\end{document}